\documentclass [a4paper,english]{article}
\usepackage{amsmath, amsfonts, amssymb}
\usepackage[cp1251]{inputenc}
\usepackage[T2A]{fontenc}
\usepackage{babel,indentfirst,epigraph}
\usepackage{graphicx,graphics,verbatim}
\newcommand{\bn}[1]{\[#1\]}
\newcommand{\be}[2]{\begin{equation}\label{#1} {#2}\end{equation}}

\newtheorem{lemme}{Lemma}
\newtheorem{coro}{Corollary}
\newtheorem{theorem}{Theorem}

\title{Proofs of some conjectures on monotonicity of  ratios of Kummer, Gauss and generalized hypergeometric functions.}
\author{ Khaled Mehrez,\; Sergei M. Sitnik }

\begin{document}
\date{}
\maketitle

\begin{abstract}
In the preprint  \cite{Sit1} one of the  authors formulated some conjectures on  monotonicity of  ratios for exponential series sections. They lead to more general  conjecture on monotonicity of  ratios of Kummer hypergeometric functions and was not proved from 1993. In this paper we prove some conjectures from   \cite{Sit1} for Kummer hypergeometric functions and its further  generalizations for Gauss and generalized hypergeometric functions. The results are also closely connected with Tur\'{a}n--type inequalities.
\end{abstract}

\section{Introduction and statement of problems.}
Let us consider the series for the exponential function
$$
\exp(x)=e^x=\sum_{k=0}^{\infty}\frac{x^k}{k!}, \ x\ge 0,
$$
its section $S_n(x)$ and series remainder $ R_n(x)$ in the form
\be{def}{S_n(x)= \sum_{k=0}^{n}\frac{x^k}{k!},
\ R_n(x)=\exp(x)-S_n(x)=\sum_{k=n+1}^{\infty}\frac{x^k}{k!}, \,x\ge 0.}

Besides simplicity and elementary nature of these functions many mathematicians studied problems for them.
G.~Szeg\H{o} proved a remarkable limit distribution for zeroes of sections, accumulated along so--called   the Szeg\H{o} curve (\cite{ESV}). S.~Ramanujan seems was the first who proved the non--trivial  inequality for exponential sections in the form (\cite{Ram}, pp. 323--324) : if
$$
\frac{e^n}{2}=R_{n-1}(n)+\frac{n^n}{n!}\theta(n)
$$
then
$$
\frac{1}{3} < \theta(n)=\frac{n!\left(\frac{e^n}{2} - R_{n-1}(n)\right)}{n^n} < \frac{1}{2}.
$$
This result is important as it also evaluates $e^n$ in rational bounds
$$
\frac{2n^n}{3n!} + 2R_{n-1}(n) < e^n <\frac{n^n}{n!} + 2R_{n-1}(n)
$$
as it was specially pointed out in (\cite{Ram}, pp. 323--324).

In the preprint  \cite{Sit1} were thoroughly studied inequalities of the form
\be{inq}{m(n)\le f_n(x)=\frac{R_{n-1}(x)R_{n+1}(x)}{\left[R_n(x)\right]^2} \le M(n), \,x\ge 0.}
The search for the best constants $m(n)=m_{best}(n),\,M(n)=M_{best}(n)$ has some history. The left--hand side of (\ref{inq}) was first proved by Kesava Menon in \cite{KM} with
$m(n)=\frac{1}{2}$ (not best) and by Horst Alzer in \cite{ALZ} with
\begin{equation}\label{mbest}
m_{best}(n)=\frac{n+1}{n+2}=f_n(0),
\end{equation}
cf. \cite{Sit1} for the more detailed history.
In  \cite{Sit1} it was also shown that in fact the inequality (\ref{inq}) with  the sharp lower constant (\ref{mbest}) is a special case of the stronger inequality proved earlier in 1982 by Walter Gautschi  in \cite{Gau1}.

It seems that the right--hand side of (\ref{inq}) was first proved by the author in \cite{Sit1} with $M_{best}=1=f_n(\infty)$. In \cite{Sit1} dozens of generalizations of inequality   (\ref{inq}) and related results were proved. May be in fact it was the first example of so called Turan--type inequality for special case of the Kummer hypergeometric functions, recently this class of inequalities became thoroughly studied (cf.  \cite{Karp1}--\cite{Lin}).

Obviously the above inequalities are consequences of the next conjecture  originally formulated in \cite{Sit1} and recently revived in \cite{Sita1}--\cite{Sita2}.

\textit{\textbf{Conjecture 1.} The function $f_n(x)$ in  (\ref{inq}) is monotone increasing\\ for $x\in [0; \infty), n\in \mathbb{N}$.} So the next inequality is valid
\be{mon}{\frac{n+1}{n+2}=f_n(0) \le f_n(x) < 1=f_n(\infty).}

In 1990's we tried to prove this conjecture in the straightforward manner  by expanding an inequality  $(f_n(x))^{'}\ge 0$ in series and multiplying  triple products of  hypergeometric functions but failed (\cite{Sit2}--\cite{Sita1}).

Consider a representation via Kummer hypergeometric functions
\be{Kum}{f_n(x)=\frac{n+1}{n+2} \ g_n(x), \ g_n(x)=\frac{_1F_1(1; n+1; x)  _1F_1(1; n+3; x)}{\left[_1F_1(1; n+2; x)\right]^2}.}
So the conjecture 1 may be reformulated in terms of this function $g_n(x)$ as conjecture 2.

\textit{\textbf{Conjecture 2.} The function $g_n(x)$ in  (\ref{Kum}) is monotone increasing\\ for $x\in [0; \infty), n\in \mathbb{N}$.}

This leads us to the next more general

\textit{\textbf{Problem 1.} Find monotonicity  in $x$ conditions for $x\in [0; \infty)$\\ for all parameters a,b,c for the function}
\be{prob}{h(a,b,c,x)=\frac{_1F_1(a; b-c; x)  _1F_1(a; b+c; x)}{\left[_1F_1(a; b; x)\right]^2}.}

We may also call (\ref{prob}) mockingly (in Ramanujan way, remember his mock theta--functions!) "\textit{The abc--problem}" for Kummer hypergeometric functions, why not?

Another generalization is to change Kummer  hypergeometric functions to higher ones.

\textit{\textbf{Problem 2.} Find monotonicity  in $x$ conditions for $x\in [0; \infty)$ for all\\ vector--valued parameters a,b,c for the function}
\be{probpq}{h_{p,q}(a,b,c,x)=\frac{_pF_q(a; b-c; x)  _pF_q(a; b+c; x)}{\left[_pF_q(a; b; x)\right]^2},}
\bn{a=(a_1,\ldots, a_p), b=(b_1,\ldots, b_q), c=(c_1, \ldots, c_q).}

This is "\textit{The abc--problem}" for generalized hypergeometric functions. The more complicated problems are obvious and may be considered for pairs or triplets of parameters and also for multivariable hypergeometric functions.

The aim of this paper is to prove conjectures 1 and 2, and  to find conditions for validity of problems 1 and 2 and so completely solve them, cf. also \cite{MS1}.

\section{Two lemmas}

We formulate two useful lemmas which will be used below. These lemmas were first proved in (\cite{BK}), cf. also (\cite{PV})--(\cite{AV}) for the detailed proof and further applications.

\begin{lemme}\label{l1}
Let $(a_{n})$ and $(b_{n})$ $(n=0,1,2...)$ be real numbers, such that $b_{n}>0,\;n=0,1,2,...$ and $\left(\frac{a_{n}}{b_{n}}\right)_{n\geq 0}$ is increasing (decreasing), then $\left(\frac{a_{0+}...+a_{n}}{b_{0}+...+b_{n}}\right)_{n}$ is increasing (decreasing).
\end{lemme}

\begin{lemme}\label{l2}.
Let $(a_{n})$ and $(b_{n})$ $(n=0,1,2...)$ be real numbers and
let the power series $A(x)=\sum_{n=0}^{\infty}a_{n}x^{n}$ and $B(x)=\sum_{n=0}^{\infty}b_{n}x^{n}$
be convergent if $|x|<r$. If $b_{n}>0,\, n=0,1,2,...$ and if the
sequence $\left(\frac{a_{n}}{b_{n}}\right)_{n\geq0}$is (strictly)
increasing (decreasing) , then the function $\frac{A(x)}{B(x)}$ is also
(strictly) increasing on $[0,r[$.
\end{lemme}

\section{Monotonicity for the Kummer hypergeometric function and associated Tur\'an type inequality }
\begin{theorem} Let $a,b,c$ be  real numbers such that $0<a<b-c$ and $b>1$ and the function $x\longmapsto h(a,b,c,x)$  is defined by
\begin{equation}
h(a,b,c,x)=\frac{_{1}F_{1}(a;b-c;x){}_{1}F_{1}(a;b+c;x)}
{\left[_{1}F_{1}(a;b;x)\right]^{2}}.
\end{equation}
Then this function  is increasing on $[0,\infty[.$ Consequently, for $n\in\mathbb{N},$ the functions $x\longmapsto f_{n}(x)$ in (\ref{inq}) and $x\longmapsto g_{n}(x)$ in (\ref{Kum}) are also increasing on $[0,\infty[.$
\end{theorem}
\textbf{Proof}
For all $a,b,c$ be  real numbers such that $0<a<b-c$ and $b>1$ we evualate
\begin{equation*}
\begin{split}
h(a,b,c,x)&=\frac{_{1}F_{1}(a;b-c;x){}_{1}F_{1}(a;b+c;x)}{\left[_{1}F_{1}(a;b;x)\right]^{2}}=\\
&=\frac{\left(\sum_{n=0}^{\infty}\frac{(a)_{n}}{(b-c)_{n}n!}x^{n}\right)
\left(\sum_{n=0}^{\infty}\frac{(a)_{n}}{(b+c)_{n}n!}x^{n}\right)}
{\left[\sum_{n=0}^{\infty}\frac{(a)_{n}}{(b)_{n}n!}x^{n}\right]^{2}}=\\
&=\frac{\sum_{n=0}^{\infty}A_{n}x^{n}}{\sum_{n=0}^{\infty}B_{n}x^{n}},
\end{split}
\end{equation*}
where
\[
A_{n}=\sum_{k=0}^{n}\frac{(a)_{k}(a)_{n-k}}{(b-c)_{k}(b+c)_{n-k}k!(n-k)!} \ \ \textrm{and}\ B_{n}=\sum_{k=0}^{n}\frac{(a)_{k}(a)_{n-k}}{(b)_{k}(b)_{n-k}k!(n-k)!}.\]
Let define sequences $(u_{n,k})_{k\geq0}$, $(v_{n,k})_{k\geq0}$ and $(w_{n,k})_{k\geq0}$  by
$$
u_{n,k}=\frac{(a)_{k}(a)_{n-k}}{(b-c)_{k}(b+c)_{n-k}k!(n-k)!},\,\,\, v_{n,k}=\frac{(a)_{k}(a)_{n-k}}{(b)_{k}(b)_{n-k}k!(n-k)!},\,\,
$$
and
$$ w_{n,k}=\frac{u_{n,k}}{v_{n,k}}=\frac{(b)_{k}(b)_{n-k}}{(b-c)_{k}(b+c)_{n-k}},\, k\geq0.
$$
It follows that
\begin{equation*}
\begin{split}
\frac{w_{n,k+1}}{w_{n,k}}&=\frac{u_{n,k+1}v_{n,k}}{v_{n,k+1}u_{n,k}}=\\
&=\frac{(b)_{k+1}(b)_{n-k-1}(b-c)_{k}(b+c)_{n-k}}{(b-c)_{k+1}(b+c)_{n-k-1}(b)_{k}(b)_{n-k}}=\\
&=\frac{\Gamma(b+k+1)}{\Gamma(b+k)}.\frac{\Gamma(b+n-k-1)}{\Gamma(b+n-k)}.\frac{\Gamma(b-c+k)}
{\Gamma(b-c+k+1)}.\frac{\Gamma(b+c+n-k)}{\Gamma(b+c+n-k-1)}=\\
&=\frac{(b+k)}{(b-c+k)}.\frac{(b+c+n-k-1)}{(b+n-k-1)}\geq1.
\end{split}
\end{equation*}

We conclude that the sequence $(w_{n,k})_{k\geq0}$ is increasing and
consequently the sequence $(C_{n}=\frac{A_{n}}{B_{n}})_{n\geq0}$ is also increasing by lemma \ref{l1}.
Thus the function $h(a,b,c,x)$ is increasing on $[0,\infty[$ by
lemma \ref{l2}. Finally, replacing $a$ and $c$ by $1$ and $b$ by $n+1$ for all $ n\in\mathbb{N}$, we obtain that the functions $x\longmapsto g_{n}(x)$ and $x\longmapsto f_{n}(x)$ are also increasing on $[0,\infty[.$ So both conjectures $1$ and $2$  from introduction are proved. And also we found the solution to the Problem $1$ from introduction if restrictions of the theorem 1 are valid.

\begin{coro} For all $a,b,c$ be  real numbers such that $0<a<b-c$ and $b>1$, the following Tur\'an type inequality
\begin{equation}\label{1ss}
\left[_{1}F_{1}(a,b,x)\right]^{2}\leq{}_{1}F_{1}(a,b-c,x).{}_{1}F_{1}(a,b+c,x)
\end{equation}
holds for all $x\in[0,\infty[.$
\end{coro}
\textbf{Proof} Since the function $x\longmapsto h(a,b,c,x)$ is increasing on $[0,\infty[$, we have
$$ h(a,b,c,x)\geq h(a,b,c,0)=1.$$

This result is interesting as a corollary of monotonicity property we consider, this inequality itself is not new and may be found in (\cite{Bar2}). And in general Tur\'an type inequalities always can be generalized to stronger results on monotonicity of function ratios with unit upper or lower constants.

\section{Monotonicity for the hypergeometric function and associated Tur\'an type inequality }

Now we also solve the Problem 2 for  general hypergeometric--type functions under some natural conditions.

\begin{theorem}
Let $p,q\in\mathbb{N}$ be such that $p\leq q+1,$ $a=(a_{1},...,a_{p}),\, b=(b_{1},...,b_{q}),\\
c=(c_{1},...,c_{q}),\, b_{i}>0,\, b_{i}-c_{i}>0$
for $i=1,2,...,q$ and $a_{i}>b_{i}$ for $i=2,...,p$. If $b_{i}>1$ for $i=1,2,...,q,$ then the function
$x\longmapsto h_{p,q}(a,b,c,x)$ in (\ref{probpq}) is strictly increasing on $[0,1[.$
\end{theorem}

\textbf{Proof}

By using the power--series representations of the function $_{p}F_{q}(a;b;x)$ we have

\begin{equation*}
\begin{split}
h_{p,q}(a;b;c;x)&=\frac{_{p}F_{q}(a;b-c;x).\,_{p}F_{q}(a;b+c;x)}{\left(_{p}F_{q}(a;b;x)\right)^{2}}=\\
&=\frac{
\left[
\sum_{n=0}^{\infty}
\frac{
(a_{1})_{n}(a_{2})_{n}...(a_{p})_{n} x^n}
{(b_{1}-c_{1})_{n}(b_{2}-c_{2})_{n}...(b_{q}-c_{q})_{n}n!}
\right]
}
{
\left[\sum_{n=0}^{\infty}\frac{(a_{1})_{n}(a_{2})_{n}...(a_{p})_{n}x^{n}
}
{(b_{1})_{n}(b_{2})_{n}...(b_{q})_{n}n!}
\right]^2
}\cdot\\
& \cdot \left[\sum_{n=0}^{\infty}\frac{(a_{1})_{n}(a_{2})_{n}...(a_{p})_{n}x^{n}}
{(b_{1}+c_{1})_{n}(b_{2}+c_{2})_{n}...(b_{q}+c_{q})_{n}n!}\right]
=\frac{\sum_{n=0}^{\infty}A_{n}(a,b,c)x^{n}}{\sum_{n=0}^{\infty}B_{n}(a,b)x^{n}}
\end{split}
\end{equation*}
where
$$
A_{n}(a,b,c)=\sum_{k=0}^{n}U_{k}(a,b,c)=
$$
$$
=\sum_{k=0}^{n}\frac{\left[(a_{1})_{k}(a_{1})_{n-k}\right]\left[(a_{1})_{k}(a_{1})_{n-k}\right]...
\left[(a_{p})_{k}(a_{p})_{n-k}\right]}{\left[(b_{1}-c_{1})_{k}...(b_{q}-c_{q})_{k}\right]
\left[(b_{1}+c_{1})_{n-k}...(b_{q}+c_{q})_{n-k}\right]k!(n-k)!},
$$
and
$$
B_{n}(a,b)=\sum_{k=0}^{n}V_{k}(a,b)=
$$
$$
=\sum_{k=0}^{n}\frac{\left[(a_{1})_{k}(a_{1})_{n-k}\right]\left[(a_{2})_{k}(a_{2})_{n-k}\right]...
\left[(a_{p})_{k}(a_{p})_{n-k}\right]}{\left[(b_{1})_{k}(b_{1})_{n-k}\right]
\left[(b_{2})_{k}(b_{2})_{n-k}\right]...\left[(b_{q})_{k}(b_{q})_{n-k}\right]k!(n-k)!}.
$$
Now, for fixed $n\in\mathbb{N}$  we define  sequences  $\left(W_{n,k}(a,b,c)\right)_{k\geq0}$ by
\begin{equation*}
W_{n,k}(a,b,c)=\frac{U_{k}(a,b,c)}{V_{k}(a,b)}=\frac{\left[(b_{1})_{k}(b_{1})_{n-k}\right]
\left[(b_{2})_{k}(b_{2})_{n-k}\right]...\left[(b_{q})_{k}(b_{q})_{n-k}\right]}
{\left[(b_{1}-c_{1})_{k}...(b_{q}-c_{q})_{k}\right]\left[(b_{1}+c_{1})_{n-k}...(b_{q}+c_{q})_{n-k}\right]}.
\end{equation*}
For $n,k\in\mathbb{N}$ we evaluate
$$
\frac{W_{n,k+1}(a,b,c)}{W_{n,k}(a,b,c)}
=\prod_{j=1}^{q}\left[\frac{(b_{j})_{k+1}(b_{j})_{n-k-1}(b_{j}-c_{j})_{k}(b_{j}+c_{j})_{n-k}}
{(b_{j})_{k}(b_{j})_{n-k}(b_{j}-c_{j})_{k+1}(b_{j}+c_{j})_{n-k-1}}\right]=
$$
$$
=\prod_{j=1}^{q}\left[\left(\frac{\Gamma(b_{j}+k+1)}{\Gamma(b_{j}+k)}\right)
\left(\frac{\Gamma(b_{j}+n-k-1)}{\Gamma(b_{j}+n-k)}\right)\right.\cdot
$$
$$
\cdot\left.
\left(\frac{\Gamma(b_{j}-c_{j}+k)}{\Gamma(b_{j}-c_{j}+k+1)}\right)
\left(\frac{\Gamma(b_{j}+c_{j}+n-k-1)}{\Gamma(b_{j}+c_{j}+n-k)}\right)\right]=
$$
$$
=\prod_{j=1}^{q}\left[\frac{b_{j}+k}{b_{j}-c_{j}+k}\right]
\left[\frac{b_{j}+c_{j}+n-k-1}{b_{j}+n-k-1}\right]>1.
$$
And now we conclude that $\left(W_{n,k}\right)_{k\geq0}$ is increasing and consequently $\left(C_{n}=\frac{A_{n}}{B_{n}}\right)_{n\geq0}$ is increasing too by the Lemma \ref{l1}. Thus the function $x\longmapsto h_{p,q}(a;b;c;x)$ is increasing on $[0,1[$ by the Lemma \ref{l2}.
It completes the proof of the theorem 2.

\begin{coro}Let $p,q\in\mathbb{N}$ be such that $p\leq q+1,$ $a=(a_{1},...,a_{p}),\, b=(b_{1},...,b_{q}),\\
c=(c_{1},...,c_{q}),\, b_{i}>0,\, b_{i}-c_{i}>0$
for $i=1,2,...,q$ and $a_{i}>b_{i}$ for $i=2,...,p$. If $b_{i}>1$ for $i=1,2,...,q,$ then the following Tur\'an type inequality
\begin{equation}\label{Fpq}
_{p}F_{q}(a;b-c;x){}_{p}F_{q}(a;b+c;x)>\left(_{p}F_{q}(a;b;x)\right)^{2}
\end{equation}
holds for all $x\in[0,1[.$
\end{coro}
\textbf{Proof}
Follows immediately from the monotonicity of the function $h_{p,q}(a;b;c;x).$

This Tur\'an type inequality (\ref{Fpq})
is included as a corollary of monotonicity property we consider, this inequality itself is not new and may be found in (\cite{Karp3}).

There are applications of considered inequalities in the theory of transmutation operators for estimating transmutation kernels and norms (\cite{Sit3}--\cite{Sit4}) and for problems of function expansions by systems of integer shifts of Gaussians  (\cite{Sit5}--\cite{Sit6}).

Recently the authors also proved generalizations of above monotonicity properties and Tur\'an type inequalities for the case of $q$--hypergeometric functions \cite{MS3}--\cite{MS4}.

The authors are thankful to D.~Karp for useful discussions.

\bigskip

Khaled Mehrez.\\
D\'{e}partement de Math\'{e}matiques,\\
l’Institut Preparatoire Aux Etudes d’Ingenieur de Monastir (IPEIM).\\
Monastir 5000, Tunisia.\\
E-mail address: k.mehrez@yahoo.fr\\

\medskip

Sergei M. Sitnik.\\
Voronezh Institute of the Russian Ministry of Internal Affairs.\\
Voronezh, Russia.\\
E-mail address: pochtaname@gmail.com
\end{document}